*Time out of mind: "Subben's checklist" revisited – A partial description of the development of quantitative OR papers over a period of 25 years*

Torbjörn Larsson
Professor of optimization, Department of Mathematics, Linköping University
torbjorn.larsson@liu.se

Michael Patriksson[1]
Professor of applied mathematics, Department of Mathematical Sciences, Chalmers University of Technology and the University of Gothenburg
mipat@chalmers.se

Johanna Pejlare
Senior lecturer, Department of Mathematical Sciences, Chalmers University of Technology and the University of Gothenburg
johanna.pejlare@gu.se

Abstract

This short paper aims to investigate some of the historical developments of *one* classic, well-cited and highly esteemed scientific journal in the domain of quantitative operations research – namely the INFORMS journal *Operations Research* – over a period of 25 years between 1981 and 2006. As such this paper, and the journal in question, represents *one* representative attempt to analyze – for the purpose of possible future generalization – how research production has evolved, and evolves, over time. Among the general developments that we think we can trace are that (a) the historical overviews (i.e., literature surveys) in the articles, as well as the list of references, somewhat counter-intuitively *shrink* over time, while (b) the motivating and modelling parts *grow*. We also attempt to characterize – in some detail – the appearance and character, over time, of the *most* cited, as well as the *least* cited, papers over the years studied. In particular, we find that many of the least cited papers are quite imbalanced; for example, some of them include one main section only, and least cited papers also have shorter reference lists.

We also analyse the articles' utilization of important "buzz words" representing the constitutive parts of an OR journal paper, based on "Subben's checklist" (Larsson & Patriksson, 2014; 2016). Based on a word count of these buzz words we conclude through a citation study, utilizing a collection of particularly highly or little cited papers, that there is a quite strong positive correlation between a journal paper being highly cited and its degree of utilization of this checklist.

## 1. Introduction – Tangled Up In Blue

In 2014, the Nordic operations research popular magazine *ORbit* published a short paper (entitled "*Subben's checklist and the quality of articles in OR*") by the first two authors of the current paper, devoted to the assessment of the characterization – and the quality – of a scientific paper in (quantitative) operations research (for short: OR). That paper was later slightly expanded, and eventually published in 2016, in the scientific journal *Computers & Operations Research*, then bearing the title *"Subben's checklist" and the assessment of articles in mathematical optimization/operations research: in memoriam of Subhash C. Narula*. Professor Subhash Narula (fondly nicknamed "Subben"), previously a professor at Rensellaer Polytechnic Institute, Troy, New York, was for a few years in the 1990s the leader of the research group in mathematical optimization at the Department of mathematics at Linköping University, Linköping, Sweden. The checklist bears his name because he was the one contributing the most items in the list, which was constructed during a long conversation in 1993 with the first two authors of this paper. The conversation, one outcome of which became the checklist, centred around possible ways in which to characterize the quality of journal articles, as well as their "completeness."

---
[1] Corresponding author

The contents of the *ORbit* and *C&OR* papers are potentially quite useful as an introductory guide for scientific authors, especially for young researchers who have not yet grasped what characterizes a scientific paper, and what makes a paper good – or very good, perhaps incomplete, or in fact plain bad!

Subben's checklist – describing the necessary parts of a complete scientific article – is as follows:

*Table 1: Subben's checklist*

| 1. Relevance | Motivation of study; need, benefit; why interesting? |
|---|---|
| 2. Background | History, references, state of the art; framework, delimitations |
| 3. Motivation | Shortcomings in existing knowledge or methodology |
| 4. Remedy | Proposal of actions in order to remove the lack of existing knowledge or methodology |
| 5. Hypothesis | Description of the research question(s) considered |
| 6. Method(ology) | Choice of research methodology |
| 7. Realization | Presentation of the new scientific contributions |
| 8. Analysis | Validation of results, conclusions, consequences, and applications; future research opportunities |

Subben's checklist, as well as the characterization of OR papers, may also fruitfully be utilized to assess important aspects of a paper, and for that it has become a quite useful tool. It has, for example, been utilized in a research methodology course in a PhD program at the Department of Technology Management and Economics at Chalmers University of Technology in Gothenburg, Sweden. It is also utilized at Mathematical Sciences, Chalmers, as a writing guide for PhD students in the area of operations research.

What has – we think – not yet been explored to a large extent is the possibility to utilize a writing guide, such as Subben's checklist, in the investigation of how the contents of articles in OR journals change over the years, and to assess – and perhaps quantify – the quality of articles. By means of one example scientific journal, namely the journal *Operations Research*, we utilize Subben's checklist in particular to illustrate how the portions of the various constitutive parts of a scientific paper, as listed in Table 1, have changed over a period of 25 years. Through the use of the checklist we believe we may trace important changes in a given journal's style and focus – and indeed general changes in quantitative operations research and its output.

In particular we investigate and compare OR papers that either are among the most cited, or among the least cited, in order to try to characterize the two "extremes". As we shall see, well-cited papers are typically more "balanced" in their contents, while the least cited papers may – for example – include essentially one section only.

The next section describes the experiment. Section 3 provides an analysis of the appearance of the constitutive parts – as listed in Table 1 –, as measured in portions of the whole, over the years studied. The analysis in particular focuses on, and contrasts, the difference between papers with either many or very few citations, in order to produce – if possible – some characteristic differences between the most well-cited and the least cited papers. Section 4 investigates the development of self-citations, and the final Section 5 provides an analysis of the results, as well as concluding remarks including potential avenues of future study.

## 2. *The study object – Things Have Changed*

First, we collected the ten *most* cited – as well as the ten *least* cited – papers in the OR journal *Operations Research* during the years 1981 and 2006, respectively (i.e., volumes 29 and 54). The journal was selected based on several criteria, among which are its consistent high quality, its long history, and the fact that it is so well-known. A time span of 25 years was deemed enough to make it possible to trace major advances in the field, changes in publication practice and editorial principles, as well as the expansion of science in general. The year 2006 was elected mainly such that enough citing years could be included up until today, for the purpose of selecting papers to bring into the study – the idea for this paper was also born in the year 2016, and the collecting started then as well.

These in total 40 papers were then scrutinized, in an attempt to characterize their content *portions*, as distinguished among the following eight categories:

- Introduction (motivation, scope),
- Review (of the literature),
- Modelling (of the problem(s) at hand),
- Theory & mathematical analysis (of the problem(s) at hand),
- Applications & practice,
- Numerical analysis,
- Future research possibilities, and
- References.

As Theory and Analysis on the one hand, Review and Future research, and Applications & practice and Numerical analysis on the other, are three quite coupled units – which may or may not be united under one banner –, they formed three pairs of single items, resulting in a final collection of six categories. (The fact that section names in the journal are not uniformly phrased also made it natural to strive for the construction of a uniform presentation.)

Table 2 is a representation of the results of a thorough scan of the 10 *most* cited papers from 1981 and 2006, respectively, and the 10 *least* cited papers from those same two years, using Scopus. The proportions of the above-mentioned parts of the papers were measured through the utilization of the word count tool in Adobe Acrobat Reader.

*Table 2: Content proportions*

| Content categories | Proportion (%) 1981 (most cited) | Proportion (%) 1981 (least cited) | Proportion (%) 2006 (most cited) | Proportion (%) 2006 (least cited) |
|---|---|---|---|---|
| Introduction, motivation | 11.0 | 13.4 | 9.5 | 14.7 |
| Survey/Review/Future research | 20.0 | 1.5 | 9.5 | 1.1 |
| Modelling | 8.0 | 42.0 | 10.0 | 14.7 |
| Theory & Analysis | 29.0 | 15.3 | 24.0 | 25.2 |
| Application, practice, and numerical analysis | 17.0 | 23.0 | 40.0 | 39.0 |
| References | 15.0 | 4.8 | 7.0 | 5.3 |

## *3. Analysis, I – The Times They Are A-Changin'*

According to the above table the 25 years from 1981 to 2006 reveal the following about the appearance of a *well-cited* paper in the journal *Operations Research*, based on our sample:

- All three categories [*Introduction, motivation*; *Survey/Review/Future research*, & *References*] associated with the *history* of the field, and *connections* to other subjects, shrink quite a lot from 1981 to 2006: collectively this portion drops almost by half, from 48% to about 26% of the total. Scientific publishing grows fast, whence there is an increasing body of references to build upon – and connect to, when creating and writing science. It is therefore quite unfortunate that there is a *diminishing connection* to the past. (As will be seen below, this category still is much more well represented compared to the case of the least cited articles.)

- The *Modelling* and *Application, practice, and numerical analysis* parts collectively grow from 25% in 1981 to 50% in 2006. Partly we think it is due to the fact that modelling "exercises" over the years have become more serious (perhaps partly as a consequence of an editorial decision?), and more importantly, realistic (in contrast to "principle models" and academic "examples"), whence more details are not only available, but in fact necessary to include in order to sufficiently well describe the applications.

- The *Theory & Analysis* section does not vary a great deal in size over the years; it remains about ¼ of the total.

- In 1981, a *Future research* section is present in five articles out of the ten, while there are in total 24 future research suggestions. In 2006 only one article omits this section; in total, there are 31 research ideas mentioned among the ten articles analysed in that year.

- The average number of *References* (already mentioned in the first item) in an article drops slightly from 53 to 45 – a decline that is perhaps not a dramatic one, while certainly over the years the volume of available literature grows.

- The average number of *words* in the articles grow over time: In 1981, the shortest paper has 7216 words, the mean is 10,176, and the longest paper has 23,252 words. In 2006, the shortest paper has 9143 words, the mean is 12,196 words, and the longest paper has 16,836 words.

For the *least cited* papers we see the following development from 1981 to 2006:

- The sections on *Introduction* and *Motivation* are slightly longer than those in the well-cited papers, and they also grow slightly over time.

- The sections on *Survey/Review* and *Future research* are always very short – no more than 1.5% of the total article content. In 1981 we find two articles, each of which presents one future research question. In 2006, we find future research ideas in six of the articles, and in total 14 research ideas (among which six stems from one single article).

- The *Theory & Analysis* sections grow from 15% to 25.2% during the 25 years.

- The *Modelling* and *Application, practice, and numerical analysis* sections collectively defines the majority of the volume of the papers: 64% in 1981, and 52.8% in the year 2006.

- The *References* section is almost always shorter than in the most well-cited papers, while it – as the average portion of the paper – even drops from 7% to 5.3% over the time period studied. In the ten least cited papers from 1981 there are 87 references, while in 2006 they are 256. (Compared with the average values of 53 and 45, respectively, in 1981 and 2006 in a well-cited paper, the least cited papers hence have less than 9 and 26 references on average, respectively, during those two years.) One must remark, however, that the font used for the reference list is smaller than in 1981, and the page itself is larger in 2006.

- The average number of *words* develops as follows: In 1981, the shortest paper has 1320 words, the mean is 5142, and the longest paper has 8622 words. In 2006, the shortest paper has 7269 words, the mean is 10,901, and the longest paper has 17,368 words. The development is hence similar to that of the most cited ones, in that there is a clear increase in volume over time. Just as the case is with the reference section, the fact that the journal style has changed between the year 1981 and 2006 – allowing more words per page – obviously affects the outcome.

From the repeated browsing and reading of the *least* cited articles, there is also – to the naked eye – a feeling that many of them not only have a narrow focus, but also are structured such that they appear to be rather imbalanced and even bordering on being *incomplete*, with – for example – the majority of some papers consisting essentially of one or two very long sections, while other common sections are either non-existing, or very short.

For reference, a cursory look at all the articles in the first issue of Operations Research in 2016 (35, respectively 10, years after the years 1981 and 2006) shows a large increase in the providing of *motivations*, and particularly an increase in material on *modelling* aspects. On the other hand, the connection to the *related literature*, *application work*, and *numerical analysis* drop dramatically, and the same can be said about the *future research* section, which is always quite short. All the while, the papers also get longer with time, although – as we see above – the reference list actually shrinks. (Further, in 2006 and 2016 the font used in the reference list is smaller, while also the actual pages are larger, as compared to the year 1981 – thus in fact allowing more information to be included on the page.)

We have not investigated the possibility that some of the changes over the years – particularly the reduction of the survey part, and the increased portion of the modelling part – are results of editorial board decisions, or an organic development based, for example, on the needs to provide more information on increasingly complex models. That would however be an interesting future study.

## 4. Analysis, II: Self-citations – I and I

We recorded self-citations (to one or more of the authors of each article) over the issues studied. (Several years ago, citation studies tended to exclude self-citations, but nowadays they typically are included in the Journal Impact Factor measure.) In the year 1981, the number of self-citations for the ten *most* well-cited papers is 91 (9.1 on average per article), and in 2006 they are 100 (10 on average per article). The analysis of the ten *least* cited papers in 1981 and 2006 reveals that there are 15 such self-citations in 1981 (1.5 on average per article), and 33 (3.3 on average per article) in 2006. (The 2016 issue has 79 self-citations – 4.8 on average per article.) In any case, from the sources utilized we *cannot* conclude that a bad habit of self-citation is present in the material used.

## 5. Analysis, III: Changing of The Guards

As an additional source of information, we searched all of the above 40 analysed articles for the number of mentions of significant words associated with the development of an article (and utilized in the paper on "Subben's checklist" – see also the table in Section 1); for ease of reference, the year with the most "hits" for a given significant word is provided in bold font:

| *Type* | min *1981* | max *1981* | min *2006* | max *2006* | Sum: |
|---|---|---|---|---|---|
| Relevance | 0 | 0 | 0 | **1** | 1 |
| Motivation | 2 | **7** | 6 | **7** | 22 |
| Shortcomings | 0 | 0 | 0 | 0 | 0 |
| Remedy | **1** | **1** | 0 | 0 | 2 |
| Hypothesis | 2 | **3** | 1 | 2 | 8 |
| Realization | 0 | 4 | 0 | **28** | 32 |
| Methodology | 2 | **16** | 7 | 13 | 38 |
| Theory | 6 | **69** | 17 | 35 | 127 |
| Validation | 0 | 0 | 0 | **28** | 28 |
| Background | 0 | 2 | 1 | **3** | 6 |
| Question | 8 | **50** | 9 | 22 | 89 |
| Consequences | 4 | 14 | 3 | **22** | 43 |
| Analysis | 25 | 11 | 21 | **64** | 121 |
| **Sum**: | 50 | **177** | 72 | **225** | |

While the "buzz words" *theory* and *questions* – words that are associated with the core of the subject, and in particular the motivation for an article's presence – dominate in 1981 (followed by *methods* and *conclusions*), in 2006 *analysis, realization* and *theory* dominate. A look at the 2016 issue reveals that there is a stronger focus on *validation* and *analysis* – words that are associated more with the results of the research, and its post-evaluation, and in 2016 *analysis* dominates even more strongly.

As an overview of the above results we remark that the *well-cited* papers over the two years studied have 402 word "hits" in Subben's checklist, while the *least cited* papers only have 112 such "hits". It does therefore appear that Subben's checklist fulfils its intended purpose quite well.

In a follow-up analysis, we complemented the search for *realization* by a search for the related term *implement*, bearing in mind the possible alternative meanings of the word. In the ten most *well-cited* papers in 1981, the word *implement* was used 30 times in the context of "algorithms", while it was used 34 times in the context of "practical realization". In the ten most *well-cited* papers in 2006, the word *implement* was used 14 times in the context of "algorithms", and never in the context of "practical realization," or "decision-making."

In the ten *least cited* papers in 1981 five papers mention the word *implement*. Among these the word relates to algorithms 16 times (13 of those reside in one article), while two articles refer to practice, in total six times. In the *least cited* papers in 2006 the word is used in eight out of the ten papers; in six articles, it is devoted to the algorithmic context only, the word-count being 25, and referring to practice in three, the word-count being ten.

In the 2016 issue *implement* refers to algorithms 16 times, and in decision-making contexts 21 times.

Can the journal Operations Research therefore be said to have become more "mathematicised" over the years? The overview in Section 3 certainly hints in that direction, partly because of an increased focus on mathematical modelling.

## 6. Evaluating science – *What Good Am I?*

The second list in the above-mentioned papers by Larsson & Patriksson (2014, 2016) refers to criteria for evaluating science. The list covers questions, research, and results from the scientific work done. It may concern all, or a subset of, the items *relevance, scientific foundation, generality, consistency, availability, scientific height and depth, originality, news value, integration, consequences, realization,* and *durability*.

While there is an abundance of material to study in order to assess these 13 criteria, a thorough study of the topic needs to be very well prepared, and it must therefore be relegated to potential future work.

## 7. Controversies in science – Ain't Talkin'

There certainly are scientific fields in which (typically senior) scientists have made strong remarks against the then current developments in scientific production, in particular regarding the item *relevance*. Among the fields of study where the authors of this paper have observed the harshest remarks are in the field of *transportation science*, where Gordon Newell (1925–2001) – a pioneer in transportation science in general, and queueing theory in particular – have stated that he favoured quality over quantity, and found that we (scientists) have failed to understand and model the behaviour of queues, and we have failed to treat travellers as they should be treated in our models, not as "a consumer good that can be sold to the highest bidder." (This and the following quotes are taken from Gordon Newell's article *Memoirs on highway traffic flow theory in the 1950s*, published in 2002. Gordon Newell feels that we have failed in developing new special techniques to the special problems that we are indeed facing, and we instead simply "rework and refine old procedures." He believes that the reason for the field not developing as strongly anymore – according to him – is that the scientists who contributed in the beginning were brilliant scholars in neighbouring fields who could bring in fresh ideas. Newell also states: "On the surface, it would seem that the subject has continued to grow and develop but actually, in my view, progress peaked in the 1960s and took a sudden dive in the 1970s." He further states: "The journal *Transportation Science* has degenerated into a journal of computer algorithms and 'optimization' relative to ad hoc objectives. There is seldom a paper dealing with some transportation issue or the answer to some question."

All of the above statements are well in line with the excellent statement by the mathematical optimization pioneer Arthur M. Geoffrion (1976):

> *The purpose of mathematical programming is insight, not numbers.*

This statement is indeed correct, and should be emphasized and discussed more often, not only among scholars, but particularly with students.

(The above statement was in fact borrowed from a quote by R.W. Hamming, 1962, pp. vii, 276, and 395, in which the phrase "mathematical programming" is replaced by "computing".) In this day and age when it appears that so much is published that it is almost impossible to know whether every journal paper includes anything actually new, these quotes are (still) timely roll calls for all conscientious researchers, reviewers, and editors to keep repetition and mediocrity at bay. A particularly timely statement that supports Hamming's and Geoffrion's is Kenneth Sörensen's paper "Metaheuristics – the metaphor exposed" (2013), which argues that the research done in that particular field threatens to become unscientific – if it isn't already.

It is therefore quite promising that the analysis made above in Section 4 at least indicates a reduction of self-promotion.

## 8. Analysis of references – "Not I", says the referee

We have investigated the appearance of the reference lists, by looking at the mean age of the references over the three years. The basis for this analysis is the investigation of whether newer papers tend to also cite relatively newer papers – that is, if the history may be "shrinking."

The ten *most* cited papers from 1981 have – in total – 530 references (that is, 53 on average), while their ages sum up to 4590 years, from their year of publication. Hence the mean age of a reference among the papers analysed in the year 1981 is 8.66 years. Performing this analysis for the ten most cited papers during the year 2006 yields the result that there are – in total – 445 references among the ten papers (that is, 44.5 on average), and the mean age is 10.8 years. The conclusion from this sample is that history is *not* shrinking.

The ten *least* cited papers from 1981 have – in total – 85 references (that is, 8.5 on average), while their ages sum up to 654 years, from their year of publication. The mean age hence is 7.7 years. The ten least cited papers in 2006 have – in total – 256 references (that is, 25.6 on average); the mean age is 12.3 years. Also in this case we see that history is *not* shrinking.

In the first issue 2016 comprising 18 articles, there are 741 references (41.2 on average), and their mean age is 12.35 years.

In this set of data, the average age of the reference list increases somewhat (by 3.7 years from 1981 to 2016). This somewhat narrow analysis may indicate that the science presented in the journal Operations Research is of a less time-dependent variety, in the sense that it does not depend on quick changes in some technology, for example. The authors of this paper certainly value the fact that history is not entirely forgotten, and in fact almost all articles cite at least some classic book or article.

In order to form some comparison with at least one other similar scientific journal in the field, we also selected 40 papers from the journal Mathematical Programming, Series A. As with the case of the journal Operations Research, we chose the ten most, and the ten least, cited papers from the years 1981 and 2006.

The ten most cited papers in Mathematical Programming, Series A in 1981 have on average 19 references, and their mean age is 8.5 years. The ten most cited papers in 2006 have on average 28 references, and their mean age is 12.7 years. The ten least cited papers in 1981 have on average 7 references, and their mean age is 9.6 years. The ten least cited papers in 2006 have on average 20 references, and their mean age is 16.75 years.

As is the case with the journal Operations Research, for the set of data utilized in the journal Mathematical Programming, Series A, we find that highly cited papers have *more* references – while the less cited ones have slightly *older* references. Why that is the case should be investigated.

## 9. *Subben's checklist and the quest for citations – Forgetful Heart*

We hypothesized in Section 1 that Subben's checklist can be utilized to measure whether a paper is well-cited. In this section, we provide an investigation on whether the number of citations to a published article is positively correlated with its degree of utilization of the checklist. A quick test from volume 54 (2006), issue 3, of Operations Research indicated that we need to incorporate synonyms to the list. Hence the below original list of eight words was appended with the following 16 terms (the original terms are given before the dash), noting that the words do have the right meaning, when scanning the papers:

Relevance – pertinence
Background – history
Motivation – incentive, reason
Remedy – solution
Hypothesis – theory, proposition, conjecture
Method, methodology – approach, technique, plan, mechanism, design, system
Realization – implementation
Analysis – assessment

As a test bed, we have again investigated the ten *most* cited papers from the years 1981 and 2006, and the ten *least* cited papers from those years. We searched those papers to see how many times the eight words in Subben's checklist are mentioned; we then also allowed the above synonyms to those words.

For the ten *most* cited papers in 1981 (the one with the least citations having 17 ditto, the highest score being 412, and the mean value being 180) we found that the number of words mentioned in Subben's checklist (and their synonyms) in those papers ranged from 3 to 14 (with a mean of 8.9), and the total number of mentions of those words were 996. The mean length of these papers is 24.7 pages, and the number of references in those papers range from 30 to 159 – with a mean value of 53.

Among the ten *least* cited papers in 1981 (three having *no* citations, and no-one more than two) we found that the number of words mentioned in Subben's checklist ranged from 2 to 7 (with a mean of 3.7), and the total number of mentions of those words were 195. The mean length of these papers is 13 pages, and the number of references in those papers range from 1 to 15, with a mean value of 9.

Among the ten *most* cited papers in 2006 (citations ranging from 46 to 216, the mean being 99) we found that the number of words mentioned in Subben's checklist ranged from 1 to 11 (with a mean of 7), and the total number of mentions of those words were 726. (The word "system" is especially popular.) The mean length of these papers is 15 pages, and the number of references in those papers range from 39 to 58, with a mean value of 45.

Among the ten *least* cited papers in 2006 (citations ranging from 2 to 16) we found that the number of words mentioned in Subben's checklist ranged from 1 to 7 (with a mean of 4), and the total number of mentions of those words were 429. The mean length of these papers is 13 pages, and the number of references in those papers range from 14 to 43, with a mean value of 26.

As a conclusion of this test, we notice that the *most* cited papers in 1981 are quite a lot longer than the ten *least* cited ones (on average 24.7 pages versus 13), while a wider range of items in Subben's checklist is also represented much more in the *most* cited papers (8.9 words per paper among the most cited, versus 3.7 for the *least* cited). In 2006 the mean length is more even between the most and least cited papers, except for the reference list which – again – is nearly double in the well-cited papers, compared to the least cited ones.

We also studied the mean age of the reference lists; the motive was to see if we can trace a move towards newer references (perhaps then illustrating swifter shifts in the subject), or not. The mean age of the references in the ten *most* cited papers in 1981 is 8.7 years, while the mean age of the *least* cited papers in 1981 is 7.7 years. The mean age of the references in the ten *most* cited papers in 2006 is 10.8 years, while the mean age of the *least* cited papers in 2006 is 12.3 years.

Hence, there appears to be no swift change in the subject matter. (A comparison with the journal Mathematical Programming, Series A, over the same years shows that the trend is similar, the latter journal having, on average, a reference list that is about two years older.)

A larger – and perhaps more detailed – study is further needed before we may state that Subben's checklist is sufficient as a tool to make a fair comparison between papers; for now, we encourage scientists to utilize Subben's checklist when preparing manuscripts, and to be generous with citations to the pertinent literature.

## 10. *Final theme: Conclusions, questions, and potential future research avenues – There's Nothing That I Wouldn't Do*

While there are books – such as the excellent one by N. Higham (1998) – that offer writing guides for the mathematical sciences, as well as good – and bad – examples of such writing, our approach and analysis concern *post-analyses* of the writing experience. In particular, we characterize – as well as we are able, utilizing the output through *two* scientific journals – the appearance of scientific papers in relation to their citation records, viewing the writing experience several years after publication and how this writing has been assessed by the readers.

Among the questions that we have asked ourselves are the following ones, to possibly be considered in future research:

1. Can the fewer references in later papers (see Section 8) partly be explained by papers more often citing surveys – rather than providing references to original papers? This question has not yet been analysed. (While general surveys do make connections to the past, they are never as detailed and precise as a dedicated one on the subject of the paper; an editor might however find it a positive feature that a reduction of survey material could yield more room for new material.)

2. Are authors more often these days citing the journals they publish in? (Some journals' editors recommend, or even *request*, that the authors try to locate pertinent articles in the same publication. Such a conduct is, however, obviously unethical.) In our limited material, we found the following: For the year 1981 – and for the ten most well-cited papers – there were 63 cites to *Oper. Res.* (among which 37 were from one paper only), while among the ten least cited papers there were five citations to the journal. In 2006, among the ten most well-cited papers there were 51 citations to *Oper. Res.*, while there were 27 citations to *Oper. Res.* among the ten least cited ones.

3. Has the style of writing changed over the years, independently of the development of the subject, and such that it may have an effect on our analysis – perhaps as a consequence of the acceleration of research and publication?

4. Is any style change explained by the fact that the subject is more developed – *matured*, as they say?

5. Who are we writing for these days – as compared to who the readers were, say, 20 years ago? For example, do scientists write in order to convey knowledge, or mainly to document their work and increase their track record? And can it be detected?

6. What, among the results of our analysis, may be transferred also to other journals? In other words, can we establish similar trends also for several other journals of a similar nature?

We hope to be able to answer at least some of these questions in due course.